\documentclass[11pt]{amsart}
\usepackage{amsfonts}
\usepackage{amsmath}
\usepackage{amsthm}
\usepackage{amssymb}
\usepackage{latexsym}
\usepackage{multicol}
\usepackage{verbatim}

\advance\textwidth by 1.2in \advance\oddsidemargin by -.6in
\advance\evensidemargin by -.6in
\newtheorem*{cor}{Corollary}
\newtheorem*{lem}{Lemma}
\newtheorem*{prop}{Proposition}

\theoremstyle{definition}
\newtheorem*{defn}{Definition}
\theoremstyle{definition}
\newtheorem*{thm}{Theorem}

\newenvironment{pf}{\proof}{\endproof}
\newcounter{cnt}
\newenvironment{enumerit}{\begin{list}{{\hfill\rm(\roman{cnt})\hfill}}{%
\settowidth{\labelwidth}{{\rm(iv)}}\leftmargin=\labelwidth%
\advance\leftmargin by
\labelsep\rightmargin=0pt\usecounter{cnt}}}{\end{list}}

\theoremstyle{remark}

\numberwithin{equation}{section} 

\def\wt{{\rm wt}}

\def\opl_#1^#2{\text{\tiny$\bigoplus\limits_{\text{\footnotesize$#1$}}^{\text{\footnotesize$#2$}}$}}

\begin{document}

\newcommand{\thmref}[1]{Theorem~\ref{#1}}
\newcommand{\secref}[1]{Section~\ref{#1}}
\newcommand{\lemref}[1]{Lemma~\ref{#1}}
\newcommand{\propref}[1]{Proposition~\ref{#1}}
\newcommand{\corref}[1]{Corollary~\ref{#1}}
\newcommand{\remref}[1]{Remark~\ref{#1}}
\newcommand{\defref}[1]{Definition~\ref{#1}}
\newcommand{\er}[1]{(\ref{#1})}
\newcommand{\id}{\operatorname{id}}
\newcommand{\tensor}{\otimes}
\newcommand{\nc}{\newcommand}
\newcommand{\rnc}{\renewcommand}
\newcommand{\qbinom}[2]{\genfrac[]{0pt}0{#1}{#2}}
\nc{\cal}{\mathcal} \nc{\goth}{\mathfrak} \rnc{\bold}{\mathbf}
\renewcommand{\frak}{\mathfrak}
\newcommand{\desc}{\operatorname{desc}}
\newcommand{\Maj}{\operatorname{Maj}}
\renewcommand{\Bbb}{\mathbb}
\nc\bomega{{\mbox{\boldmath $\omega$}}} \nc\bpsi{{\mbox{\boldmath
$\Psi$}}}
 \nc\balpha{{\mbox{\boldmath $\alpha$}}}
 \nc\bpi{{\mbox{\boldmath $\pi$}}}

\nc{\krm}{KR(m\omega_1)}

\newcommand{\lie}[1]{\mathfrak{#1}}
\makeatletter
\def\section{\def\@secnumfont{\mdseries}\@startsection{section}{1}%
  \z@{.7\linespacing\@plus\linespacing}{.5\linespacing}%
  {\normalfont\scshape\centering}}
\def\subsection{\def\@secnumfont{\bfseries}\@startsection{subsection}{2}%
  {\parindent}{.5\linespacing\@plus.7\linespacing}{-.5em}%
  {\normalfont\bfseries}}
\makeatother
\def\subl#1{\subsection{}\label{#1}}
 \nc{\Hom}{\text{Hom}}
\nc{\ch}{\text{ch}} \nc{\ev}{\text{ev}}

\nc{\krsm}{KR^\sigma(m\omega_2)}
\nc{\krsmzero}{KR^\sigma(m_0\omega_i)}
\nc{\krsmone}{KR^\sigma(m_1\omega_i)}
 \nc{\vsim}{v^\sigma_{2,m}}
 \nc{\gr}{{\rm gr}}

 \nc{\Cal}{\cal} \nc{\Xp}[1]{X^+(#1)} \nc{\Xm}[1]{X^-(#1)}
\nc{\on}{\operatorname} \nc{\Z}{{\bold Z}} \nc{\J}{{\cal J}}
\nc{\C}{{\bold C}} \nc{\Q}{{\bold Q}}
\renewcommand{\P}{{\cal P}}
\nc{\N}{{\Bbb N}} \nc\boa{\bold a} \nc\bob{\bold b} \nc\boc{\bold
c} \nc\bod{\bold d} \nc\boe{\bold e} \nc\bof{\bold f}
\nc\bog{\bold g} \nc\boh{\bold h} \nc\boi{\bold i} \nc\boj{\bold
j} \nc\bok{\bold k} \nc\bol{\bold l} \nc\bom{\bold m}
\nc\bon{\bold n} \nc\boo{\bold o} \nc\bop{\bold p} \nc\boq{\bold
q} \nc\bor{\bold r} \nc\bos{\bold s} \nc\bou{\bold u}
\nc\bov{\bold v} \nc\bow{\bold w} \nc\boz{\bold z} \nc\boy{\bold
y} \nc\ba{\bold A} \nc\bb{\bold B} \nc\bc{\bold C} \nc\bd{\bold D}
\nc\be{\bold E} \nc\bg{\bold G} \nc\bh{\bold H} \nc\bi{\bold I}
\nc\bj{\bold J} \nc\bk{\bold K} \nc\bl{\bold L} \nc\bm{\bold M}
\nc\bn{\bold N} \nc\bo{\bold O} \nc\bp{\bold P} \nc\bq{\bold Q}
\nc\br{\bold R} \nc\bs{\bold S} \nc\bt{\bold T} \nc\bu{\bold U}
\nc\bv{\bold V} \nc\bw{\bold W} \nc\bz{\bold Z} \nc\bx{\bold x}

\nc{\krmi}{KR(m\omega_i)} \nc{\krsmi}{KR^\sigma(m\omega_i)}
\nc{\vim}{v_{i,m}} \nc{\vsimi}{v^\sigma_{i,m}} {\title[]{
Kirillov--Reshetikhin modules associated to  $G_2$}
\author{Vyjayanthi Chari and  Adriano Moura}
\address{Department of Mathematics, University of
California, Riverside, CA 92521.} \email{chari@math.ucr.edu}\thanks{VC was partially
supported by the NSF grant DMS-0500751}
\address{UNICAMP - IMECC, Campinas - SP - Brazil, 13083-970.} \email{aamoura@ime.unicamp.br}

\maketitle

\begin{abstract}
We define and study the Kirillov--
Reshetikhin
modules for algebras of type $G_2$. We compute the graded
character of these modules and verify that they are in 
accordance
with the conjectures in \cite{HKOTY}, \cite{HKOTT}. These 
results
give the first complete description of families of
Kirillov--Reshetikhin modules whose isotypical components 
have
multiplicity bigger than one.
\end{abstract}

\section*{Introduction}

In \cite{cmkr} we defined and studied a family of restricted
modules for the current and twisted current algebras associated to
a finite--dimensional classical simple Lie algebra $\lie g$ and a
diagram automorphism of $\lie g$ of order two. These modules, which we
called the restricted Kirillov-Reshetikhin modules,  are given by
generators and relations and were denoted by $KR^\sigma(m\omega_i)$, where $\sigma$ is the diagram automorphism,
$i$ is a node  of the Dynkin diagram of the subalgebra $\lie g_0$ of $\lie g$ consisting of the fixed points of $\sigma$, and $m$ is a non--negative integer.
Here we understand $\sigma$ to be the identity in the untwisted case.
They admit a natural grading which is compatible with the grading on the current algebras. In
particular, the graded pieces are finite--dimensional  modules for $\lie g_0$.
It was proved in \cite{cmkr} that, regarded as $\lie g_0$-modules,  there were
no non--zero maps between the distinct graded pieces  and, moreover,
 the multiplicity of any irreducible representation
 in a particular graded piece was at most one.
In fact, the graded character was computed in \cite{cmkr} and
verified to be in accordance with the conjectures in \cite[Appendix
A]{HKOTY} and \cite[Section 6]{HKOTT} for the usual Kirillov-Reshetikhin modules for the corresponding quantum affine algebras.
 When $\lie g_0$ is an exceptional Lie algebra,
the conjectures in these papers  make it clear that for some nodes
of the Dynkin diagram one or both of the aforementioned properties of
the graded pieces may fail. The modules $KR^\sigma(m\omega_i)$ are
known to be isomorphic  to the Demazure modules, further details
can be found in \cite{CL}, \cite{cmkr}, \cite{FoL1}, \cite{FoL2}.

In this paper we define and study the modules for the current
algebra associated to $G_2$ and to the twisted current algebra
associated to $D_4$ and a diagram automorphism of order three. In
both cases the fixed point subalgebra   $\lie g_0$ is of type $G_2$.
We prove that the conjectures of \cite{HKOTY} and \cite{HKOTT} are
true in these cases. In particular, there are now maps of $\lie
g_0$--modules between the distinct non--zero graded pieces for
$KR^\sigma(m\omega_i)$ for some $i$ and the multiplicity of an
irreducible module in a graded piece can be greater than  one.
Moreover,  our result on the graded character of the module
$KR(m\omega_1)$ for $G_2$ is actually an improvement on the
conjectural graded--character formula in \cite{HKOTY} which has
some multiplicity--zero  terms.

The overall scheme of the proof is very similar to the one in
\cite{cmkr}: we prove that the conjectural character formula is an
upper bound for the character and then we prove that it is also a
lower bound. However, one runs into difficulty almost immediately
as the underlying combinatorics is rather more complicated. In
order to prove the upper bound we use an elementary but useful
result on representations of the 3--dimensional Heisenberg
algebra. For the lower bound, as in \cite{cmkr}, we first study
some ``fundamental'' Kirillov-Reshetikhin modules and then
realize the other modules as a submodule of a tensor product of
the fundamental ones. But this time the fundamental modules are
too big to be constructed explicitly as in \cite{cmkr}. To solve
this we use the notion of fusion product of modules of the current
algebra, which was introduced and studied in \cite{FL},
\cite{FKL}. The second step, in  which involves studying graded
quotients of tensor products of the fundamental
Kirillov--Reshetikhin modules, is really much more complicated,
since one has to prove not only that a particular representation
occurs in a given grade, but also one has to determine its
multiplicity. Identifying these quotients  and proving that the
isotypical components occur is non--trivial, since the projection
of the natural vectors do not generate the desired $\lie
g_0$--submodule. To solve this part
  we use the explicit description of some
highest--weight vectors in tensor products of representations of
$\lie{sl}_2$ and in tensor products of fundamental representations
for $\lie g_0$.

The paper is organized as follows. In section \ref{prel} we fix
the basic notation and collect the results we will need for the
proof. In section \ref{definition} we define the
Kirillov-Reshetikhin modules, state the main theorem, and make the
connection with the conjectures in \cite{HKOTY} and \cite{HKOTT}.
We prove the theorem in sections \ref{ubs} and \ref{lbs}.

\setcounter{section}{0}
\section{Preliminaries}\label{prel}

\subsection{The Lie algebra $G_2$ and its representations}\label{notsimple}
Throughout this paper $\lie
g_0$ will denote the Lie algebra of type $G_2$, $\lie h_0$ a
Cartan subalgebra of $\lie g_0$ and  $\alpha_i$, $i=1,2$, a set of
simple roots where we assume that $\alpha_1$ is short and
$\alpha_2$ is long.
Let $R_l^+$ and $R_s^+$ be the set of positive long and positive
short roots respectively,
$$R^+_l=\{\alpha_2,\alpha_2+3\alpha_1, 2\alpha_2+3\alpha_1\},\ \
R_s^+=\{\alpha_1,\alpha_1+\alpha_2, 2\alpha_1+\alpha_2\}.$$ Given
$\alpha\in R^+$ we denote by $x^\pm_\alpha$ any non-zero element
of $(\lie g_0)_{\pm \alpha}$.  The subalgebras $\lie n_0^\pm$ are
defined in the obvious way by  $\lie n_0^\pm=\oplus_{\alpha\in
R^+}\bc x^\pm_\alpha.$  Let $\omega_i$, $i=1,2$, be the
fundamental weights and note that $\omega_1=2\alpha_1+\alpha_2$
and $\omega_2=2\alpha_1+3\alpha_2$. Let $P$ (resp. $Q$) be the
integer lattice spanned by the fundamental weights (resp. simple
roots) and let $P^+$ (resp. $Q^+$) be the $\bz_+$ span of the
fundamental weights (resp. simple roots). Fix elements
$h_{\alpha_i}$, $i=1,2$, such that
$\omega_j(h_{\alpha_i})=\delta_{ij}$ for $i,j=1,2$. Then it is
easy to see that $[x^+_{\alpha_i}, x^-_{\alpha_i}]$ is a non--zero
multiple of $h_{\alpha_i}$.

Given a finite--dimensional $\lie g_0$--module $V$, we have
$$V=\oplus_{\lambda\in P}^{}V_\lambda,\ \ V_\lambda=\{v\in V:
hv=\lambda(h)v\ \forall\ h\in\lie h\}.$$ Let $\wt(V)=\{\mu\in P:
V_\mu\ne 0\}$ and given $0\ne v\in V_\mu$ set $\wt(v)=\mu$. Let $Z[P]$
be the integral group ring of $P$ with basis $e(\mu)$, $\mu\in P$,
and set
$$\ch(V)=\sum_{\mu\in P}\dim(V_\mu)e(\mu).$$For
$\lambda\in P^+$, let $V(\lambda)$ be the irreducible $\lie
g_0$--module with highest weight $\lambda$ and highest weight
vector $v_\lambda$. Thus $V=\bu(\lie g_0)v$, where $$\lie
n_0^+v_\lambda=0, \ \ hv=\lambda(h)v,\ \
(x_{\alpha_i}^-)^{\lambda(h_i)+1}v =0.$$ Note that
$$
\ch(V(\omega_1))=e(0)+\sum_{\pm\alpha\in R^+_s}e(\alpha),\ \
\ch(V(\omega_2))=2e(0)+\sum_{\alpha\pm\in R^+}e(\alpha).$$

We shall need the following result which is trivially proved.
\begin{lem}\label{hwvg2} Given $0\le p\le s\in\bz_+$, there exists $a,b\in\bc^\times $ such that the following holds in
$V(\omega_2)^{\otimes s}\otimes V(\omega_1)$:
$$\lie n_0^+\left((x^-_{\alpha_1+\alpha_2}v_2^{\otimes s-p})\otimes v_2^{\otimes p}\otimes v_1+
\left( ax^-_{\alpha_1+\alpha_2}+b
x_{\alpha_2}^-x_{\alpha_1}^-\right) (v_2^{\otimes s}\otimes
v_1)\right)=0.$$ \hfill\qedsymbol

\end{lem}

\subsection{The associated current  and twisted current algebras} Given a Lie algebra $\lie a$, let  $\lie a[t]=\lie
a\otimes \bc[t]$ be the polynomial current algebra of $\lie a$
with bracket $[x\otimes t^r, y\otimes t^s]=[x,y]\otimes t^s$.

From now on, let  $\lie g$ be a Lie algebra  of type $D_4$ and
$\sigma$  the automorphism of $\lie g$ induced by an automorphism
of order three of the Dynkin diagram. Let $\xi$ be a primitive
cube root of unity. Then,
$$\lie g=\opl_{j=0}^2\lie g_j,\ \ \lie g_j=\{x\in\lie g:
\sigma(x)=\xi^j x\}.$$ Notice that the notation $\lie g_0$ is
unambiguous since it is known that the fixed point subalgebra of
$\sigma$ is isomorphic to $G_2$. Further, the subspaces $\lie
g_r$, $r=1,2$, are clearly representations of $\lie g_0$ and in
fact $\lie g_r\cong V(\omega_1)$. For $\alpha\in R^+_s$, we let $y_\alpha^\pm$, $z^\pm_\alpha$ be
{\bf{any}} non--zero elements in $(\lie g_1)_{\pm \alpha}$ and
$(\lie g_2)_{\pm \alpha}$, respectively.

Extend $\sigma$ to an automorphism $\sigma_t$ of $\lie g[t]$ by
$x\otimes t^s\to \sigma(x)\otimes \xi^st^s$. Let $\lie
g[t]^\sigma$ be the set of fixed points of $\sigma_t$. Then,
$$\lie g[t]^\sigma =\lie g_0\otimes\bc[t^3]\oplus \lie g_1\otimes t\bc[t^3]\oplus\lie g_2\otimes t^2\bc[t^3].$$
Set $$\lie n^\pm[t]^\sigma=\lie n_0^\pm\otimes
\bc[t^3]\opl_{\alpha\in R_s^+}^{} \left(\bc y^\pm_{\alpha}\otimes
t\bc[t^3]\oplus \bc z^\pm_{\alpha}\otimes t^2\bc[t^3]\right).$$

 We shall use the fact that as
vector spaces
$$\bu(\lie g_0[t])\cong \bu(\lie n_0^-[t])\bu(\lie h_0[t])\bu(\lie
n_0^+[t]),\ \ \bu(\lie g[t]^\sigma)\cong \bu(\lie
n^-[t]^\sigma)\bu(\lie h[t]^\sigma)\bu(\lie n^+[t]^\sigma)$$
without further comment.

\subsection{Graded modules and graded characters} The algebras
$\lie g_0[t]$ and $\lie g[t]^\sigma$ are obviously $\bz_+$--graded
algebras. Given a $\bz_+$--graded module $V_t=\oplus_{n\in\bz_+}
V_t[n]$ for $\lie g_0[t]$ or $\lie g[t]^\sigma$, it is easy to see
that $V_t[n]$ is a $\lie g_0$--module.  If $V_t[n]$ is
finite--dimensional for all $n\in\bz_+$, the graded character of
$V_t$ is defined by
$$\ch_t(V_t)=\sum_{n\in\bz_+}t^n\ch(V_t[n])=\sum_{n\in\bz_+}t^n\left (\sum_{\mu\in
P^+}m_{\mu,n}(V_t[n])\ch(V(\mu))\right),$$
where $m_{\mu,n}(V_t[n])$ are given by
 $$V_t[n]\cong_{\lie g_0}\opl_{\mu\in P^+}^{} V(\mu)^{\oplus m_{\mu,n}(V_t[n])}.$$

Set
$$V_t(n)=V_t/\opl_{s>n}^{}V_t[s].$$

We end this section with some results  which are  used crucially
later in the paper.

\subsection{A result on representations of $\lie sl_2$} Let $x^+,x^-,h$ be
 the standard basis for the Lie algebra $\lie sl_2$ and let
  $V(s)$ be the $(s+1)$--dimensional representation of $\lie{sl}_2$
with highest weight vector $v_s$.

\begin{lem}\label{sl2}
Given $0\le p\le s\in\bz_+$ and
$j\le \min(p,s-p)$, there exist $c_1,\cdots ,c_j\in\bc^\times $
such that the following holds in $V(1)^{\otimes s}$.
$$x^+\left(\left({(x^-)}^j v_1^{\otimes s-p}\right)\otimes v_1^{\otimes p}+ \sum_{\ell=1}^{j} c_\ell(x^-)^\ell
\left(\left((x^-)^{j-\ell}v_1^{\otimes s-p}\right)\otimes v_1^{\otimes
p}\right)\right)=0.$$ \end{lem}
\begin{pf} Notice first that the $\lie sl_2$--module of $V(1)^{\otimes s}$
generated by $v^{\otimes s-p}_1\otimes v_1^{\otimes p}$ can be
identified with  the  submodule of $V(s-p)\otimes V(p)$ generated
by $v_{s-p}\otimes v_p$. Hence it suffice to prove that for $0\le
j\le \min(p,s-p)$ there exist $c_1,\cdots ,c_j\in\bc^\times  $
such that
$$x^+\left(\left({(x^-)}^j v_{s-p}\right)\otimes v_p+ \sum_{\ell=1}^{j} c_\ell(x^-)^\ell
\left(\left((x^-)^{j-\ell}v_{ s-p}\right)\otimes v_{
p}\right)\right)=0.$$ But this is immediate from the
Clebsch--Gordon formulas.
\end{pf}

\subsection{A  result on representations of the Heisenberg  algebra}

\begin{lem}\label{heis}  Consider the three dimensional Heisenberg algebra $\frak
H$ spanned by elements $x,y,z$ where $z$ is central and $[x,y]=z$.
Suppose that $V$ is a  representation of $\frak H$ and let $0\ne
v\in V$ be such that $x^rv=0$. Then  for all $k,s\in\bz_+$  the
element $y^kz^sv$ is in the span of elements of the form
$x^ay^bz^cv$ with $a>0$ and $0\le c<r$.\end{lem}
\begin{pf} Suppose first that $r=1$, then $$y^kz^sv=
y^kz^{s-1}xyv= y^kxyz^{s-1}v =xy^{k+1}z^{s-1}v-ky^{k}z^sv,
$$
i.e.,
$$ (k+1)y^{k}z^sv = xy^{k+1}z^{s-1}v.$$

Assume that we know the result for $r'<r$. Then, we have
$$y^kz^sv=
y^kz^{s-1}xyv- y^kz^{s-1}yxv= xy^{k+1}z^{s-1}v-ky^{k}z^sv -
y^{k+1}z^{s-1}xv,$$
i.e.,
$$(k+1)y^kz^sv=xy^{k+1}z^{s-1}v-y^{k+1}z^{s-1}xv.$$

 Since $x^{r-1}xv=0$ it
follows by induction that $y^{k+1}z^{s-1}xv$ is in the span of
elements of the form $x^ay^bz^cxv$ with $a>0$ and $c<r-1$. But such
elements are clearly in the span of elements of the form
$x^ay^bz^cv$ with $a>0$ and $c<r$. An induction on $s$ again gives
the result.
\end{pf}

\subsection{Fusion Products} We shall need the following result
which was proved in \cite{FL}, \cite{FKL}. We state it in
the form and in the case of interest to us.

\begin{prop}\label{fus} Let $V_i$, $i=1,2$, be finite--dimensional  graded $\lie g_0[t] $ modules generated by elements
$v_i$, $i=1,2$, satisfying the relations: $\lie n_0^+v_i=0$,
$h\otimes t^r=\delta_{r,0}\lambda_i(h)v_i$ for some $\lambda_i\in
P^+$, $i=1,2$, and  all $h\in\lie h_0$. Then, there exists a graded
$\lie g_0[t]$--module denoted $V_1*V_2$ which is generated by an
element $v$ satisfying:
$$\lie n_0^+ v=0,\ \ h\otimes t^r=\delta_{r,0}(\lambda_1+\lambda_2)(h)v,
\ \ \ \forall\ h\in \lie h_0$$ and
$$V_1*V_2\cong_{\lie g_0} V_1\otimes V_2.$$
\hfill\qedsymbol
\end{prop}

\section{The Kirillov--Reshetikhin modules}\label{definition}

In this section we define and prove some elementary properties of
the Kirillov--Reshetikhin modules for the algebras of type $\lie
g_0[t]$ and $\lie g[t]^\sigma$.

\subsection{The KR--modules for $\lie g_0[t]$}
\begin{defn} For  $m\in\bz_+$, let $KR(m\omega_i)$ be
the $\lie g_0[t]$--module generated by an element $v_{i,m}$ with
relations,
\begin{equation}\label{deftwu}
\lie n_0^+[t]v_{i,m}=0,\ \ (h\otimes
t^s)v_{i,m}=\delta_{s,0}m\omega_i(h)v_{i,m},\ \
\end{equation}
for all $h\in\lie h_0$, $s\in\bz_+$, and
\begin{equation}\label{deftw2} x^-_{\alpha_j}v_{i,m}=0,\ \ i\ne
j,\ \  (x^-_{\alpha_i})^{m+1}v_{i,m}=0,\ \ (x_{\alpha_i}^-\otimes
t)v_{i,m}=0.
\end{equation}
\hfill\qedsymbol\end{defn}

\subsection{The KR--modules for $\lie g[t]^\sigma$}

\begin{defn} For  $m\in\bz_+$, let $KR^\sigma(m\omega_i)$ be
the $\lie g[t]^\sigma$--module generated by an element
$v^\sigma_{i,m}$ with relations,
\begin{equation}\label{deftw13}
\lie n^+[t]^\sigma v^\sigma_{i,m}=0,\ \ (h\otimes
t^s)v^\sigma_{i,m}=\delta_{s,0}m\omega_i(h)v^\sigma_{i,m},\ \
 \end{equation} for all
$h\otimes t^s \in \lie h[t]^\sigma$,and
\begin{equation}\label{deftw4}x^-_{\alpha_j}v^\sigma_{i,m}=(x^-_{\alpha_i})^{m+1}v^\sigma_{i,m}=
  (y_{\alpha_1}^-\otimes t)v^\sigma_{i,m} = (x^-_{\alpha_2}\otimes t^3)v^\sigma_{i,m} =0,
\end{equation}
where $j\ne i$.\hfill\qedsymbol
\end{defn}

\subsection{The Main Theorem} Let $\boe_i\in\bz_+^4$, $1\le i\le 4$ be the standard
basis and set $\mathbf 0=(0,0,0,0)$. Define
$\wt,\wt^\sigma:\bz_+^4\to P^+$ and $\gr,\gr^\sigma:\bz_+^4\to \bz_+$ by
\begin{align}
&\wt(\bor)=(m-r_1-3r_2-3r_3)\omega_1+(r_2+r_3-r_4)\omega_2,\ &&\gr(\bor)=r_1+r_2+2r_3+2r_4,\\
&\wt^\sigma(\bor)=(r_1+r_2-r_3)\omega_1+(m-r_1-r_2-r_4)\omega_2,\
&&\gr^\sigma(\bor)=r_1+2r_2+2r_3+3r_4,
\end{align}
where $\bor =(r_1,r_2,r_3,r_4)$.
Set
\begin{equation}
\cal A_1=\{\bor\in\bz_+^4: r_4\le r_2,\ 2r_1+3r_2+3r_3\le m\},
\end{equation}
\begin{equation}
\cal A_2^\sigma=\{\bor\in\bz_+^4: r_3\le r_1,\ \ r_1+r_2+r_3+r_4\le m\}.
\end{equation}

 The main result of this paper is the following.
\begin{thm} \label{main} Let $m\in\bz_+$. The modules $\krmi$ and $\krsmi$ are
$\bz_+$--graded and
\begin{alignat*}{3}\tag{{\it i}}
&\ch_t(KR(m\omega_1))=\sum_{\bor\in\cal A_1}t^{\gr(\bor)} \ch(V(\wt(\bor))), \ \
&&\ch_t(KR(m\omega_2))=\sum_{r=0}^m t^{m-r} \ch(V(r\omega_2)),\\ \tag{{\it ii}}
&\ch_t(KR^\sigma(m\omega_1))=\sum_{r=0}^m t^{m-r} \ch(V(r\omega_1)), \
&&\ch_t(KR^\sigma(m\omega_2))=\sum_{\bor\in\cal A_2^\sigma} t^{\gr^\sigma(\bor)} \ch(V(\wt^\sigma(\bor))).
\end{alignat*}
\end{thm}

As a consequence of the proof of the theorem we also have:

\begin{cor} \begin{enumerit}
\item Let $m\in\bz_+$, $m=3m_0+m_1$. The module $KR(m\omega_1)$ is
isomorphic to the  submodule of $ KR(3\omega_1)^{\otimes
m_0}\otimes KR(m_1\omega_1)$ generated by the element
$v_{1,3}^{\otimes m_0}\otimes v_{1,m_1}.$ \item The module
$KR(m\omega_2)$ (resp. $KR^\sigma(m\omega_i)$, $i=1,2$) is
isomorphic to the  submodule of $KR(\omega_2)^{\otimes m}$ (resp.
$(KR^\sigma(\omega_i))^{\otimes m}$) generated by the element
$v_{2,1}^{\otimes m}$ (resp.  $(v^\sigma_{i,1})^{\otimes m}$).
\end{enumerit}\hfill\qedsymbol

\end{cor}

 We prove the theorem in the next two
sections.
\subsection{The connection with the conjectures  in \cite
{HKOTY}, \cite{HKOTT}}\label{conj}

The following formulas were conjectured in \cite[Appendix A]{HKOTY} and \cite[Section 6]{HKOTT}.
\begin{align*}
&\ch_t(KR(m\omega_1)) = \sum_{k=0}^{[m/3]}\sum_{j=2k}^{m-k}
p_{j,k}(t) \ch(V((m-j-k)\omega_1+k\omega_2)),\\
&\ch_t(KR(m\omega_2)) = \sum_{k=0}^{m} t^{m-k} \ch(V(k\omega_2)),\\
&\ch_t(KR^\sigma(m\omega_1)) = \sum_{k=0}^{m} t^{m-k} \ch(V(k\omega_1)),\\
&\ch_t(KR^\sigma(m\omega_2)) = \sum_{\begin{subarray}{c}
        j,k \in \bz_+\\ j + k \le m
    \end{subarray}}
p^\sigma_{j,k}(t)\ch(V(j \omega_1+k\omega_2)),
\end{align*}
where $[s]$ denotes the biggest integer smaller than or equal to $s$,
$$p_{j,k}(t)=\left(1+[\frac{j-2k}{3}]+\min(0, [\frac{m+k-2j}{3}])\right) \sum_{s=0}^k t^{j-k+s},$$
and
$$p^\sigma_{j,k}(t)= \left(1+ \min(k,m-j-k)\right)\sum_{s=0}^{j} t^{3m-2j-3k+s}.$$

For the modules $KR(m\omega_2)$ and $KR^\sigma(m\omega_1)$ it is
clear that  Theorem \ref{main} establishes the conjectures.

In order to establish the conjecture for $KR(m\omega_1)$,  write
$m=3m_0+m_1$ with $m_1=0,1,2$. Define an equivalence relation on
$\bz^4$ by $\bor\sim\bor'$ iff $\wt(\bor)=\wt(\bor')$ and
$\gr(\bor)=\gr(\bor')$. It is easy to see that
\begin{equation*}\label{eqg2}
\bor\sim\bor' \quad\text{iff}\quad \bor-\bor' =
\ell(3\boe_1-\boe_2-\boe_4)\quad\text{for some }\ell\in\bz.
\end{equation*}
Let $\bar\bor$ be the equivalence class of $\bor$.

Let $j,k\in\bz_+$ be such that $0\le k\le m_0=[m/3]$, $2k\le j\le
m-k$. Assume also that $p_{j,k}(t)\ne 0$, i.e., $r_4+k \le
m_0+[(m_1-2r_1)/3]$, where $r_4\in\bz_+$ and $0\le r_1\le 2$ are defined by $r_4 = [(j-2k)/3]$ and $j-2k = r_1+3r_4$.
For $0\le s\le k$, set 
$$\bor_{j,k,s} = (r_1,k+r_4-s,s,r_4).$$
 It is easy to check that,
\begin{enumerit} 
\item $\bor_{j,k,s}\in\cal A_1$
\item $\wt(\bor_{j,k,s}) =(m-j-k)\omega_1+k\omega_2$, \ \ $\gr(\bor_{j,k,s})=j-k+s$, 
\item $\#\ \bar\bor_{j,k,s}\cap\cal A_1=1+[\frac{j-2k}{3}]+\min(0,[\frac{m+k-2j}{3}]).$
\end{enumerit}
Here $\#S$ is the cardinality of the set $S$. In other words, we see that
$$p_{j,k}(t)=\sum_{k=0}^{m_0}\sum_{j=2k}^{m-k}\sum_{s=0}^k(\#\ \bar\bor_{j,k,s}\cap\cal A_1) t^{\gr(\bor_{j,k,s})}.$$
Thus to show that  the conjecture in \cite{HKOTY} coincides with
Theorem \ref{main} in the case of $KR(m\omega_1)$, it suffices to
show that 
$$\{\bor_{j,k,s}: 0\le s\le k\le m_0,\ \ 2k\le j\le m-k\}$$ 
is a complete set of representatives for the equivalence classes of $\cal A_1$, i.e.,
$$\cal A_1=\bigcup_{0\le s\le k\le m_0\atop{2k\le j\le m-k}} \bar\bor_{j,k,s}\cap\cal A_1.$$
But this is now easy to do.

In the case of $KR^\sigma(m\omega_2)$ we proceed similarly.
Namely, we define an equivalence relation on $\cal \bz^4$ by
$\bor\sim\bor'$ iff $$\wt^\sigma(\bor)=\wt^\sigma(\bor'),\ \
\gr^\sigma(\bor)=\gr^\sigma(\bor')$$ and we let $\bar\bor$ be the
equivalence class of $\bor$.

It is easy to see that
\begin{equation*}
\bor\sim\bor'\iff\
\bor-\bor'=\ell(\boe_1+\boe_3-\boe_4)\quad\text{for some }\ell\in\bz.
\end{equation*}
Given $j,k,s\in\bz_+$ satisfying $j+k\le m$, $0\le s\le j$, set
$\bor_{j,k,s}=(j-s,s,0,m-j-k)$. Then,
\begin{enumerit}
\item $\bor_{j,k,s}\in\cal A_2^\sigma$
\item $\wt^\sigma(\bor_{j,k,s}) =j\omega_1+k\omega_2$, \ \ $\gr^\sigma(\bor_{j,k})=3m-2j-3k+s$, 
\item $\#\ \bar\bor_{j,k,s}\cap \cal A_2^\sigma= 1+\min(k,m-j-k).$
\end{enumerit} In other words, we see that
$$p_{j,k}(t)=\sum_{j+k=0}^{m}\sum_{s=0}^j(\#\ \bar\bor_{j,k,s}\cap \cal A_2^\sigma) t^{\gr(\bor_{j,k,s})}.$$
Thus to show that  the conjecture in \cite{HKOTY} coincides with
Theorem \ref{main} in the case of $KR^\sigma(m\omega_2)$, it
suffices to show that
$$\{\bor_{j,k,s}: 0\le j+k\le m,\ \ 0\le s\le j\},$$ 
is a complete set of representatives for
the equivalence classes of $\cal A^\sigma_2$, i.e.,
$$\cal A^\sigma_2=\bigcup_{0\le j+ k\le m\atop{0\le s\le j}} \bar\bor_{j,k,s}\cap\cal
A^\sigma_2,$$ which is easily done. 

\subsection{} We conclude this section with some elementary properties of the
modules $\krmi$ and $\krsmi$.  The proof of the next proposition
is standard (see \cite{cmkr}) and we omit the details.

\begin{prop}\label{eleprop}
Let $m\in\bz_+$, $i=1,2$, and assume that $K_m$ (resp. $v_m$) is
either  $\krmi$ or $\krsmi$ (resp. $\vim$ or $\vsimi$).

 \begin{enumerit}
\item $K_0\cong \bc$.

\item For all $\alpha\in R_0^+$, we have
$$(x^-_\alpha\otimes 1)^{m\omega_i(h_\alpha)+1} v_{m}=0.$$

\item We have
$$K_m=\opl_{\mu\in\lie h_0^*}^{} (K_m)_\mu $$ and $(K_m)_\mu\ne 0$ only if
$\mu\in m\omega_i-Q_0^+$.

\item Regarded as a $\lie g_0$--module, $K_m$ and $K_m[s]$,
$s\in\bz_+$, are isomorphic to a direct sum of irreducible
finite--dimensional representations.

\item For all $0\le r\le m$,  there exists a canonical
homomorphism $K_m\to K_r\otimes K_{m-r}$ of graded $\lie
g_0[t]$--modules (resp. $\lie g[t]^\sigma$--modules) such that
$v_m\mapsto v_{r}\otimes v_{m-r}$.

\end{enumerit}
\end{prop}

\begin{cor}$$\ch_t(K_m)=\sum_{r\in\bz_+}t^r\left(\sum_{
\mu\in P^+} m_{\mu,r}(K_m)\ch(V(\mu))\right)$$ for some
$m_{\mu,r}(K_m)\in\bz_+$.
\end{cor}

\subsection{} The next lemma is easily deduced (see \cite{cmkr}) from the defining
relations of the modules $\krmi$ and $\krsmi$.

\begin{lem}\label{red1} Let $m\in\bz_+$, $i=1,2$. Let $\alpha \in R_0^+$ and assume that
$\alpha=s_i\alpha_i+s_j\alpha_j$, $i\ne j$.
\begin{enumerit}
\item In $\krmi$ we have $$(x^-_\alpha\otimes t^r)v_{i,m}=0\ \ \forall \ \ r\ge s_i.$$

\item In $KR^\sigma(m\omega_2)$ we have
$$(x^-_\alpha\otimes t^{3r})v^\sigma_{2,m}\ =\ (y^-_\alpha\otimes t^{3r+1})v^\sigma_{2,m}\ = \
(z^-_\alpha\otimes t^{3r+2})v^\sigma_{2,m}\ =\ 0 \ \
\forall \  r\ge s_2.$$

\item In $KR^\sigma(m\omega_1)$ we have
$$(x^-_\alpha\otimes t^{3s})v^\sigma_{1,m}\ =\
(y^-_\alpha\otimes t^{3r-2})v^\sigma_{1,m}\ =\
(z^-_\alpha\otimes t^{3s-1})v^\sigma_{1,m}\ =\ 0 \ \
\forall \  r\ge s_1,\ s\ge \min(1,s_1).$$
\end{enumerit}
Here we set $y^-_\alpha=z^-_\alpha=0$ if $\alpha$ is  a long root.\hfill\qedsymbol
\end{lem}

\section{Upper bounds}\label{ubs}

\subl{}  The main result of this section is the following.

\begin{prop}\label{ub} Let $\mu\in P^+$, $k\in\bz_+$.
\begin{enumerit}
\item We have
$$m_{\mu, k}(KR(m\omega_1))\le \#\{\bor\in\cal A_1: (\mu,k)=(\wt(\bor),\gr(\bor))\},$$

$$m_{\mu,k}(KR(m\omega_2))\le 1 \quad\text{and}\quad m_{\mu,k}(KR(m\omega_2))=0 \quad\text{if}\quad \mu\ne (m-k)\omega_2.$$

\item We have
$$m_{\mu,k}(KR^\sigma(m\omega_1))\le 1 \quad\text{and}\quad m_{\mu,k}(KR^\sigma(m\omega_1))=0 \quad\text{if}\quad \mu\ne (m-k)\omega_1.$$

$$m_{\mu, k}(KR^\sigma(m\omega_2))\le \#\{\bor\in\cal A^\sigma_2: (\mu,k)=(\wt^\sigma(\bor),\gr^\sigma(\bor))\}.$$
\end{enumerit}
\end{prop}

The proposition is proved in the rest of this section.

\subsection{The case of $KR(m\omega_2)$ and $KR^\sigma
(m\omega_1)$}

 We fix  an ordered basis of $\lie n_0^-[t]$
as follows: the basis consists of elements in the set
$$\{x^-_{\alpha}\otimes t^s: \alpha\in R^+, s\in\bz_+\},$$ with
any total order that satisfies $x^-_\alpha\otimes
t^s<x^-_\beta\otimes t^r$ if $s<r$ for all $\alpha,\beta\in R^+$.
An application of the PBW theorem and Lemma \ref{red1}(i) shows
that
$$KR(m\omega_2)=\sum_{r\in\bz_+}\bu(\lie
g_0)(x^-_{3\alpha_1+2\alpha_2}\otimes t)^rv_{2,m},$$ and that
$\lie n_0^+(x^-_{3\alpha_1+2\alpha_2}\otimes t)^rv_{2,m}=0$. This
immediately implies that
$$KR(m\omega_2)=\opl_{r\in\bz_+}^{} t^rV((m-r)\omega_2)^{\oplus m_r},$$ where $0\le m_r\le 1$.

We fix  an ordered basis of $(\lie n^-)[t]^\sigma $ as follows:
the basis consists of elements in the set
$$\{X^-_{\alpha}\otimes t^{s}: \alpha\in R^+, s\in\bz_+\}, $$ where $X^-_\alpha\in\{x^-_\alpha, y^-_\alpha, z^-_\alpha\}$
and $s$ are such that $X^-_\alpha\otimes t^s\in\lie g[t]^\sigma$
with any total order that satisfies $X^-_\alpha\otimes
t^s<X^-_\beta\otimes t^r$ if $s<r$ for all $\alpha,\beta\in R^+$.
Using Lemma \ref{red1}(iii) and the Poincare Birkhoff--Witt basis
we see that
$$KR^\sigma(m\omega_1)=\sum_{r\in\bz_+}\bu(\lie
g_0)(y^-_{2\alpha_1+\alpha_2}\otimes t)^rv^\sigma_{1,m},\ \ \lie
n_0^+(y^-_{2\alpha_1+\alpha_2}\otimes t)^rv^\sigma_{1,m}=0, $$ and
the proposition  follows as before  in this case.

\subsection{The case of $KR(m\omega_1)$ and
$KR^\sigma(m\omega_2)$}  We now fix  an ordered basis of $\lie
n_0^-[t]$ as follows: the basis consists of elements in the set
$$\{x^-_{\alpha}\otimes t^s: \alpha\in R^+, s\in\bz_+\}.$$ Fix any
total  order on this set that satisfies the following:
\begin{enumerit}
\item for all $\alpha,\beta\in R^+$ and $r>0$, we have
$x^-_\alpha<x^-_\beta\otimes t^r$,
\item  further, we have

$$x^-_{3\alpha_1+2\alpha_2}\otimes t^2<
x^-_{3\alpha_1+\alpha_2}\otimes t^2<
x^-_{3\alpha_1+\alpha_2}\otimes t< x^-_{2\alpha_1+\alpha_2}\otimes
t< x^-_{3\alpha_1+2\alpha_2}\otimes t< x^-_\beta\otimes t^s,$$ for
all $(\beta,s)$ with $\beta\in R^+$, $s>0$ and $$(\beta,
s)\notin\{(3\alpha_1+2\alpha_2,2), (3\alpha_1+\alpha_2,2),
(3\alpha_1+2\alpha_2,1), (3\alpha_1+\alpha_2,1),
(2\alpha_1+\alpha_2,1)\}.$$
\end{enumerit}

 For $\lie n^-[t]^\sigma$ we adopt a similar notation.
Set   $X^-_\beta\in\{x^-_\beta, y^-_\beta,
z^-_\beta\}$,\begin{enumerit} \item for all $\alpha,\beta\in R^+$
and $r>0$, we have $x^-_\alpha<X^-_\beta\otimes t^r$, \item
further, we have
$$x^-_{3\alpha_1+2\alpha_2}\otimes
t^3<z^-_{2\alpha_1+\alpha_2}\otimes
t^2<z^-_{\alpha_1+\alpha_2}\otimes
t^2<y^-_{2\alpha_1+\alpha_2}\otimes
t<y^-_{\alpha_1+\alpha_2}\otimes t<X^-_\beta\otimes t^s,$$
for all $(\beta,s)$ with $\beta\in R^+$, $s>0$ and $$(\beta,
s)\notin\{(3\alpha_1+2\alpha_2,3), (2\alpha_1+\alpha_2,2),
(\alpha_1+\alpha_2,2), (2\alpha_1+\alpha_2,1),
(\alpha_1+\alpha_2,1)\}.$$
\end{enumerit}

  Given $\bor\in\bz_+^4$,  let
$\boy_\bor\in \bu(\lie g_0[t])$ and $\boy_\bor^\sigma\in
\bu(\lie g[t]^\sigma)$ be defined by
$$\boy_\bor=(x^-_{3\alpha_1+2\alpha_2}\otimes
t^2)^{r_4}(x^-_{3\alpha_1+\alpha_2}\otimes
t^2)^{r_3}(x^-_{3\alpha_1+\alpha_2}\otimes
t)^{r_2}(x^-_{2\alpha_1+\alpha_2}\otimes t)^{r_1}$$
and
$$\boy_\bor^\sigma=(x^-_{3\alpha_1+2\alpha_2}\otimes
t^3)^{r_4}(z^-_{2\alpha_1+\alpha_2}\otimes
t^2)^{r_3}(z^-_{\alpha_1+\alpha_2}\otimes
t^2)^{r_2}(y^-_{\alpha_1+\alpha_2}\otimes t)^{r_1},$$
respectively. If
$\bor\notin\bz_+^4$, then we set $\boy_\bor=0$ (resp.
$\boy^\sigma_\bor=0$).

 Using Lemma
\ref{red1} and the PBW theorem we see that
$$KR(m\omega_1)=\sum_{s\in\bz_+, \bor\in\bz_+^4}\bu(\lie n_0^-)
\boy_\bor(x^-_{3\alpha_1+2\alpha_2}\otimes t)^s v_{1,m},$$
$$KR^\sigma(m\omega_2)=\sum_{s\in\bz_+, \bor\in\bz_+^4}\bu(\lie
n_0^-)\boy^\sigma_\bor (y^-_{2\alpha_1+\alpha_2}\otimes
t)^sv^\sigma_{2,m}.$$ It is easy to see that the relations
$x^-_{\alpha_2}v_{1,m}=0$ and  $x^-_{\alpha_1}v^\sigma_{2,m}=0$
imply the following:
$$(x^-_{3\alpha_1+2\alpha_2}\otimes
t)^s(x^-_{3\alpha_1+\alpha_2}\otimes t)^{r_2}
v_{1,m}\in\bc\left((x^-_{\alpha_2})^s
(x^-_{3\alpha_1+\alpha_2}\otimes t)^{s+r_2}v_{1,m}\right),\ $$  $$
(y^-_{2\alpha_1+\alpha_2}\otimes
t)^s(y^-_{\alpha_1+\alpha_2}\otimes t)^{r_1}v^\sigma_{2,m}\in\bc
 \left((x^-_{\alpha_1})^s(y^-_{\alpha_1+\alpha_2}\otimes t)^{s+r_1}v^\sigma_{2,m}\right).$$
 Since
 $$[x^-_{\alpha_2},\boy_{\bor-r_2\boe_2}]\in\sum_{\bos\in\bz_+^4}\bc\boy_\bos,\
 \ \
 [x^-_{\alpha_1},\boy^\sigma_{\bor-r_1\boe_1}]\in\sum_{\bos\in\bz_+^4}\bc\boy^\sigma_\bos,$$
we see that $$\boy_\bor(x^-_{3\alpha_1+2\alpha_2}\otimes t)^s
v_{1,m}\in \sum_{\bos\in\bz_+^4}\bu(\lie n^-_0)\boy_\bos
v_{1,m},$$ $$\ \ \boy^\sigma_\bor(y^-_{2\alpha_1+\alpha_2}\otimes
t)^sv^\sigma_{2,m}\in\sum_{\bor\in\bz_+^4}\bu(\lie
n^-_0)\boy_\bos^\sigma v_{2,m}.$$ In other words we have proved
that
 \begin{eqnarray}\label{ub1}& KR(m\omega_1)&=\sum_{\bor\in\bz_+^4}\bu(\lie n^-_0)\boy_\bor v_{1,m}
 =\sum_{\bor\in\bz_+^4}\bu(\lie g_0)\boy_\bor v_{1,m}, \\
\label{ub2} &KR^\sigma(m\omega_2)&=\sum_{\bor\in\bz_+^4}\bu(\lie
n^-_0)\boy_\bor^\sigma v^\sigma_{2,m} = \sum_{\bor\in\bz_+^4}\bu(\lie
g_0)\boy_\bor^\sigma v^\sigma_{2,m} .\end{eqnarray}

\subsection{} \begin{lem}\label{nplusaction} Let $k\in\bz_+$.
\begin{enumerit}\item We have,
\begin{equation}
(x_1^+)^k\boy_\bor v_{1,m}\in\sum_{j=0}^k \bc\boy_{\bor+(k-3j)\boe_1+(j-k)\boe_2+j\boe_4} v_{1,m},\ \
(x_2^+)^k \boy_\bor v\in\bc\boy_{\bor+k\boe_3-k\boe_4}v_{1,m}.
\end{equation}
 The elements $\{\boy_\bor v_{1,m}:\bor\in\bz_+^4\}$ span a finite--dimensional
representation of $\lie n_0^+$ and hence $KR(m\omega_1)$ is a
finite--dimensional $\lie g_0[t]$--module. \item We have,
\begin{equation}\label{na}(x_1^+)^k\boy^\sigma_\bor
v^\sigma_{2,m}\in\bc
\boy^\sigma_{\bor-k\boe_3+k\boe_2}v^\sigma_{2,m},\ \ \ (x_2^+)^k
\boy^\sigma_\bor v^\sigma_{2,m}\in\sum_{j=0}^k\bc
\boy^\sigma_{\bor-(2k+j)\boe_1+(k-j)\boe_3+j\boe_4}v^\sigma_{2,m}.
\end{equation}The elements $\{\boy^\sigma_\bor v^\sigma_{2,m}:\bor\in\bz_+^4\}$
span a finite--dimensional representation of $\lie n_0^+$ and
hence $KR^\sigma(m\omega_2)$ is a finite--dimensional $\lie
g[t]^\sigma$--module.
\end{enumerit}
\end{lem}
\begin{pf} We prove (ii), the proof of (i) is identical.
The observation that
$$[x^+_{\alpha_1}, \boy^\sigma_{\bor-r_3\boe_3}]=0,\ \ [x_{\alpha_1}^+,
z_{2\alpha_1+\alpha_2}\otimes t^2]\in\bc
z_{\alpha_1+\alpha_2}\otimes t^2,$$ proves the first inclusion in
\eqref{na}. To prove the second, we begin by observing that
$$[x_2^+, z_{2\alpha_1+\alpha_2}\otimes t^2]=0$$ and
$$[x_2^+,x^-_{\theta}\otimes t^3]\in\bc x^-_{3\alpha_1+\alpha_2}\otimes t^3,\ \
[x_2^+,z_{\alpha_1+\alpha_2}\otimes t^2]\in\bc z_{\alpha_1}\otimes
t^2,\ \ [x_2^+,y_{\alpha_1+\alpha_2}\otimes t]\in\bc
y_{\alpha_1}\otimes t.$$
Lemma \ref{red1} and the commutation relations in
$\lie g[t]^\sigma$ now prove that for any $\bos\in\bz_+^4$,
$$(x^-_{3\alpha_1+\alpha_2}\otimes t^3)\boy^\sigma_\bos v_{2,m}^\sigma=0, \ \
(z_{\alpha_1}\otimes t^2) \boy^\sigma_\bos v^\sigma_{2,m} =0.$$

 Since $[x_2^+,
(y_{\alpha_1+\alpha_2}\otimes t)^r]$ is in the span of the
elements $(y_{\alpha_1+\alpha_2}\otimes
t)^{r-1}(y_{\alpha_1}\otimes t)$,\linebreak
 $(z_{2\alpha_1+\alpha_2}\otimes t^2)(y_{\alpha_1+\alpha_2}\otimes t)^{r-2}$ and
$(x^-_{\theta}\otimes t^3)(y_{\alpha_1+\alpha_2}\otimes
t)^{r-3})$, we find that \eqref{na} follows for $k=1$ from a
further application of Lemma \ref{red1}.

In particular, we have shown that
the subspace spanned by the elements $\{\boy^\sigma_\bor
v^\sigma_{2,m}:\bor\in\bz_+^4\}$ is a representation of $\lie
n_0^+$.  To see that the subspace is finite--dimensional, note
that for each $\mu\in P$, the set
$\{\bor\in\bz_+^4:\wt(\boy_\bor)=\mu\}$ is finite. Hence if the
subspace was  infinite--dimensional, there would exist an infinite
family of elements $\boy_{\bor_j}$, $j\ge 1$ with
$\wt(\boy_{\bor_j})\ne \wt(\boy_{\bor_k})$ if $j\ne k$. Since
$KR^\sigma(m\omega_2)$ is a direct sum of finite--dimensional
irreducible $\lie g_0$--modules, it follows that there must exist
an infinite family of distinct elements $\mu_j\in P^+$ such that
$KR^\sigma(m\omega_2)_{\mu_j}\ne 0$. But this is impossible since
there are only finitely many elements in the $ m\omega_2-Q^+$. The
fact that $KR^\sigma(m\omega_2)$ is finite--dimensional is
immediate from \eqref{ub1}.
\end{pf}

\subsection{}\label{ubpf} Let $\pi_0: KR(m\omega_1)\to \bu(\lie g_0)v_{1,m}$ be
the canonical  projection of $\lie g_0$--modules so that we have
$KR(m\omega_1)=\bu(\lie g_0)v_{1,m}\oplus  ker(\pi_0)$. If $\pi_0$
is injective, the proposition is proved. Otherwise there exists
$\bor_1\in\bz_+^4$ such that the element $\boy_{\bor_1}v_{1,m}$
has a non--zero projection onto $ker(\pi_0)$. Moreover, $\bor_1$
can be chosen so that: $\boy_\bos v_{1,m}\in \bu(\lie g_0)v_{1,m}$
if $\bos\in\bz_+^4$ is such that either $\wt(\bos)-\wt(\bor_1)\in
Q^+\backslash\{0\}$
or $\wt(\bos)=\wt(\bor_1)$ with $\bos<\bor_1$,
where $<$ is the lexicographic ordering on $\bz_+^4$ given by
$$(r_1,r_2,r_3,r_4)< (s_1,s_2,s_3,s_4)\iff \ \ r_{k}<s_{k},\ \text{ where }\   k=\min\{1\le p\le 4: r_p\ne s_p\}.$$
Let $v_1$ be the projection of $\boy_{\bor_1}v_{1,m}$ onto
$ker(\pi_0)$. Using Lemma \ref{nplusaction} we see that $$\lie
n_0^+ \boy_{\bor_1}v_{1,m}\in \bu(\lie g_0)v_{1,m},$$ and hence
$\lie n_0^+ v_{1,m}=0$.

Repeating this  argument, we see that we
can choose $\bor_0,\cdots,\bor_k\in\bz_+^4$ and elements $v_j\in
KR(m\omega_1)_{\wt(\bor_j)}$, $0\le j\le k$ such that:
\begin{enumerit}
\item $\wt(\bor_0)=m\omega_1\ge \wt(\bor_1)\ge\cdots\ge
\wt(\bor_k),$ \item $\lie n_0^+v_j=0$, $0\le j\le k$,
$v_0=v_{1,m}$
\end{enumerit} such that the following holds:
\begin{enumerit} \item[(a)] as $\lie g_0$--modules $KR(m\omega_1)=\oplus_{j=0}^k\bu(\lie
g_0)v_j\cong \oplus_{j=0}^k V(\wt(\bor_j))$,  \item[(b)] the
projection of $\boy_{\bor_j}$  onto $\bu(\lie g_0)v_j$ is $v_j$.
Moreover if $\bos\in\bz_+^4$ is such that either
$\wt(\bos)-\wt(\bor_j)\in Q^+$ or  $\wt(\bos)=\wt(\bor_j)$, with
$\bos<\bor_j$, then
$\boy_{\bos}v_{1,m}\in\oplus_{p=0}^{j-1}\bu(\lie g_0)v_p.$

\end{enumerit} Proposition \ref{ub} is proved for $KR(m\omega_1)$
if we show that $\bor_j\in\cal A_1$ for all $0\le j\le k$. We
first prove that if $\bor_j=(r_1,r_2,r_3,r_4)$, then $r_4\le r_2$.
For this, note that
$$x_{\alpha_2}^+(x^-_{3\alpha_1+\alpha_2}\otimes t)^{r_2}(x^-_{2\alpha_1+\alpha_2}\otimes
t)^{r_1}v_{1,m}=0,\ \
(x^-_{\alpha_2})^{r_2+1}(x^-_{3\alpha_1+\alpha_2}\otimes
t)^{r_2}(x^-_{2\alpha_1+\alpha_2}\otimes t)^{r_1}v_{1,m}=0.$$ The
subalgebra of $\lie g_0[t]$ spanned by
$(x^-_{3\alpha_1+2\alpha_2}\otimes t^2)$,
$(x^-_{3\alpha_1+\alpha_2}\otimes t^2)$ and $x^-_{\alpha_2}$ is
isomorphic to the three dimensional Heisenberg algebra. Lemma
\ref{heis}  now implies that if $r_4>r_2$, then
$\boy_{\bor_j}v_{1,m}$ is in the span of elements of the form
$(x^-_{\alpha_2})^a\boy_{\bos}v_{1,m}$ with $a>0$ and
$\wt(\bos)>\wt(\bor_j)$.  But such elements have zero projection
on $\bu(\lie g_0)v_j$ and hence $\boy_{\bor_j}$ has zero
projection onto $\bu(\lie g_0)v_j$ which contradicts (b).

Next suppose that there exists $0\le j\le k$ such that
$\bor_j=(r_1,r_2,r_3,r_4)$ and
 $2r_1+3r_2+3r_3>m$. Setting,
$\bos=\bor_{j_0}+r_1\boe_2-r_1\boe_1$, we see from Lemma
\ref{nplusaction} that
\begin{equation}\label{ls}(x_{\alpha_1}^+)^{r_1}\boy_{\bos}v_{1,m} = \boy_{\bor_{j_0} }
v_{1,m}+  \sum_{p=1}^{r_1}\boy_{\bor_{j_0}-p(3\boe_1-\boe_2-\boe_4)}v_{1,m}.
\end{equation}
Now, $\bor_{j_0}-p(3\boe_1-\boe_2-\boe_4)<\bor_{j_0}$ if $p\ge 1$
 it follows that the projection of $\boy_{\bor_{j_0}-p(3\boe_1-\boe_2-\boe_4)}v_{1,m}$ for $p\ge 1$ onto $\bu(\lie
g_0)v_j$ is zero. Since $ \boy_{\bor_{j_0} } v_{1,m}$ has a
non--zero projection onto $\bu(\lie g_0)v_{j_0}$, we see using
\eqref{ls} that $\boy_{\bos}v_{1,m}$  also has a non--zero
projection onto $\bu(\lie g_0)v_{j_0}$. Now,
$$\wt(\bos)=\wt(\bor_{j_0})-r_1\alpha_1=
(m-3r_1-3r_2-3r_3)\omega_1+(r_1+r_2+r_3-s_4)\omega_2.$$ Since
$2r_1+3r_2+3r_3>m$ it follows that $\wt(\bos)$ is not
dominant integral and so we must have that $$\wt(\bos)+
(3r_3+3r_2+3r_1-m)\alpha_1\in\wt(\bu(\lie g_0)v_{j_0})\subset
\wt(\bor_{j_0})-Q^+,$$ i.e.,
$$\wt(\bor_{j_0})+(2r_1+3r_2+3r_3-m)\alpha_1\in \wt(\bor_{j_0})-Q^+$$
which is impossible. Proposition \ref{ub} is proved for
$KR(m\omega_1)$. The result is deduced for $KR^\sigma(m\omega_2)$
in exactly the same way. One works with the Heisenberg algebra
spanned by $x_{\alpha_1}^-$,  $z^-_{\alpha_1+\alpha_2}\otimes t^2$
and $z^-_{2\alpha_1+\alpha_2}\otimes t^2$ and we omit the details.

\section{Lower Bounds}\label{lbs}

\subsection{}
The main result of this section is the
following Proposition which  together with Proposition \ref{ub} proves Theorem
\ref{main}.
\begin{prop}\label{lb}\hfill
\begin{enumerit}
 \item We have
$$m_{\mu, k}(KR(m\omega_1))\ge \#\{\bor\in\cal A_1: (\mu,k)=(\wt(\bor),\gr(\bor))\},$$

$$m_{(m-k)\omega_2, k}(KR(m\omega_2))\ge 1.$$

\item We have
$$ m_{(m-k)\omega_1, k}(KR^\sigma(m\omega_1))\ge 1.$$

$$m_{\mu, k}(KR^\sigma(m\omega_2))\ge \#\{\bor\in\cal A^\sigma_2: (\mu,k)=(\wt^\sigma(\bor),\gr^\sigma(\bor))\}.$$
\end{enumerit}
\end{prop}

\subsection{The modules $KR(m\omega_2)$}\label{om2}

Note that the $\lie g_0$ module $V(\omega_2)$ is isomorphic to the
adjoint representation of $\lie g_0$. Let $<,>$ be the Killing
form on $\lie g_0$. If $m=1$, then it is straightforward to check
that the formulas
$$(x\otimes t^r)(y,a)=(\delta_{r,0}[x,y], \delta_{r,1}<x,y>),$$
define a graded $\lie g_0[t]$--module structure on $K= V(\omega_2)
\oplus\bc$ with $K[0]=V(\omega_2)$, $K[1]=\bc$. It is trivial to
check that $K$ is a $\lie g_0[t]$ module quotient of
$KR(\omega_2)$ which proves the proposition for $m=1$. Moreover,
the assignment $x^+_{3\alpha_1+2\alpha_2} \mapsto v_{2,1}$ extends
to a $\lie g_0[t]$--module isomorphism $K\cong KR(\omega_2)$ and
hence,
$$(x^-_{3\alpha_1+2\alpha_2}\otimes t)v_{2,m}\ne 0,\ \
\lie n_0^+(x^-_{3\alpha_1+2\alpha_2}\otimes t)v_{2,m}=0.$$

 For $m>1$, consider the module
$$K_m= K(0)^{\otimes m-k}\otimes K^{\otimes k}.$$ Let $\bar
v_{2,m}$ be the image of $v_{2,m}$ in $K(0)$ and set $$\bar
v_m=\bar v_{2,m}^{\otimes m-k}\otimes v_{2,m}^{\otimes k}.$$
Using the explicit description of the module $K$, it is now easy
to see that the module $\bar K_m=\bu(\lie g_0[t])\bar v_m$ is a
quotient of $KR(m\omega_2)$. Moreover,
$$(x^-_{3\alpha_1+2\alpha_2}\otimes t)^k\bar v_m\ne 0,\ \ \lie
n_0^+((x^-_{3\alpha_1+2\alpha_2}\otimes t)^k\bar v_m)=0,$$ which
proves that $m_{(m-k)\omega_2,k}(\bar K_m)\ne 0$. Since
$KR(m\omega_2)$ is a semisimple $\lie g_0$--module it follows that
$m_{(m-k)\omega_2,k}(KR(m\omega_2))\ne 0$, thus proving the
proposition.

\subsection{The modules $KR(m\omega_1)$, $1\le m\le 3$.}

Using Proposition \ref{ub}, we see that
$m_{\mu,k}(KR(\omega_1))=0$ if $k\ne0$. Since the
formula
$$(x\otimes t^r)v=\delta_{r,0}xv,\ \ \forall\ \  x\in\lie g_0,\ \ v\in
V(\omega_1),$$ defined a $\lie g_0[t]$--module action on
$V(\omega_1)$ which makes it  a quotient of $KR(\omega_1)$, we are
done.

For $m=2$, we see from Proposition \ref{ub} that
\begin{equation}\label{2om1} m_{\mu,k}(KR(2\omega_1))=0,\ \ (\mu, k)\notin
\{ (2\omega_1,0)\ (\omega_1,1)\}.\end{equation}
Consider the fusion
product $K=KR(\omega_1)*KR(\omega_1)$. Using Proposition
\ref{fus}, we see that  $\bar v= (x^-_{\alpha_1}\otimes
t)(v_{1,m}*v_{1,m})$ is a non--zero element of $K$ and moreover,
$$\lie n_0^+[t]\bar v=0, \ \ (\lie h_0\otimes t\bc[t])\bar v=0,\ \
h\bar v=\omega_2(h)\bar v,\ \ (x^-_{\alpha_2}\otimes t)\bar v=
0,$$ for all $h\in \lie h_0$. In other words, $\bar K=\bu(\lie
g_0[t])\bar v$ is a graded  $\lie g_0[t]$--module quotient of
$KR(\omega_2)$, and it follows from Section \ref{om2} that either
$$\ch_t(\bar K)=\ch( V(\omega_2))  \ \ {\text{ or }}\  \ ch_t(\bar
K)= \ch(V(\omega_2))+t\ch(\bc).$$ Since
$$K\cong_{\lie g_0} V(2\omega_1)\oplus V(\omega_2)\oplus
V(\omega_1)\oplus \bc,$$ it follows that either $$K/\bar K\cong_{\lie
g_0} V(2\omega_1)\oplus V(\omega_1) \ \ {\text{or}} \ \ K/\bar
K\cong_{\lie g_0} V(2\omega_1)\oplus V(\omega_1)\oplus \bc. $$ An
application of Proposition \ref{fus} again shows that $K/\bar K$
is a quotient of $KR(2\omega_1)$. Equation \eqref{2om1} implies
that
$$K/\bar K\cong_{\lie g_0} V(2\omega_1)\oplus V(\omega_1),\ \
KR(2\omega_1)\cong K/\bar K,$$ and proves Proposition \ref{lb} in
this case.

For $m=3$, we see from Proposition \ref{ub} that
\begin{equation}\label{3om1} m_{\mu,r}(KR(3\omega_1))=0,\ \ (\mu, r)\notin\{
(3\omega_1,0)\ (2\omega_1,1),\ (\omega_2,1),\ (\omega_2, 2),\
(0,3) \}.\end{equation}

Consider the fusion product $K=KR(\omega_2)*KR(\omega_2)$. Set
$\bar K=\bu(\lie g_0[t])(x^-_{\alpha_2}\otimes
t)(v_{2,1}*v_{2,1})$. Using Proposition \ref{fus} one checks
easily that
 and
$\bar K$ is a quotient of $KR(3\omega_1)$ and that $K/\bar K$ is a
quotient of $KR(2\omega_2)$. Since
$$
\bar K\oplus K/\bar K\cong_{\lie g_0} K\cong_{\lie g_0}
V(2\omega_2)\oplus V(3\omega_1)\oplus V(2\omega_1)\oplus 3
V(\omega_2)\oplus 2\bc, $$ the  proposition follows for $m=3$ from
equation \ref{3om1} together with the fact that
$$K/\bar K\subset_{\lie g_0} KR(2\omega_2)\cong_{\lie g_0}
V(2\omega_2)\oplus V(\omega_2)\oplus \bc.$$

\subsection{} We shall need the following result.
\begin{lem}\label{fundvan} Let $\bor\in\bz_+^4$.
\begin{enumerit} \item In $KR(\omega_1)$ we have $\boy_\bor v_{1,1}=0$
for all $\bor\in\bz_+^4$, $\bor\ne \mathbf 0$.
\item In $KR(2\omega_1)$
we have $$\boy_\bor v_{1,2}=0\ \ \iff\ \  \bor\notin\{\mathbf 0,
\boe_1,\boe_2\},$$ and
$$\lie n_0^+ \boy_{\boe_1}v_{1,2}=0,\ \ \boy_{\boe_2}v_{1,2}\in\bc x^-_{\alpha_1}\boy_{\boe_1}v_{1,2}.$$
\item In $KR(3\omega_1)$ we have $$\boy_\bor v_{1,3}=0\ \ \iff \ \
\bor\notin\{\mathbf 0, \boe_1,\boe_2, \boe_3, \boe_4, 2\boe_1, 2\boe_2,
\boe_1+\boe_2,\boe_2+\boe_4, 3\boe_1\}.$$ The elements
$\boy_{\boe_2}v_{1,3}$, $x_{\alpha_1}^-\boy_{\boe_1}v_{1,3}$ are
linearly independent and there exists $a\in\bc^\times$, such that
$$\lie n_0^+\boy_{\boe_1}v_{1,3}=\lie
n_0^+\boy_{\boe_3}v_{1,3}=\lie n_0^+\boy_{3\boe_1}v_{1,3}=\lie
n_0^+ (\boy_{\boe_2}-ax_{\alpha_1}^-\boy_{\boe_1})v_{1,3}=0.$$
Finally,
$$\boy_{\boe_2+\boe_4}v_{1,3}\in\bc(
\boy_{3\boe_1}v_{1,3}),\ \ \boy_{\boe_4}v_{1,3} \in\bc
(x_{\alpha_2}^-\boy_{\boe_3}v_{1,3}), \ \ \boy_{2\boe_1}v_{1,3}
\in\bc (x_{\alpha_1+\alpha_2}^-\boy_{\boe_3}v_{1,3}).$$
\end{enumerit}

\end{lem}
\begin{pf} Part (i) is obvious. For (ii), it is  clear from the
fact  $ch_t(KR(2\omega_1))=\ch(V(2\omega_1))+t\ch (V(\omega_1))$
that $\boy_\bor v_{1,m}=0$ if $\bor\notin\{\mathbf 0, \boe_1,\boe_2\}$.
For the converse, suppose that $\boy_{\boe_1}v_{1,2}=0$. Since
$\wt(\boe_2)<\omega_1$, this means that if $\bor\in\bz_+^4$ is
such that $\wt(\bor)=\omega_1$, then $\boy_{\bor}v_{3,m}\in
V(2\omega_1)$  and hence proves that $m_{\mu, 1}(KR(2\omega_1))=0$
which is a contradiction. A simple calculation proves that
$x_{\alpha_1}^+\boy_{\boe_2}v_{1,2}\in\bc^{\times}(\boy_{\boe_1}v_{1,2})$
and hence it follows that $\boy_{\boe_2}v_{1,2}\ne 0$. The second
equality in part (ii) is trivially established. The proof of (iii)
is a similar detailed analysis based on the graded character of
$\ch_t(KR(3\omega_1))$.
\end{pf}

\subsection{The modules $KR(m\omega_1)$, $m>3$}

Set $K = KR(3\omega_1)$, $v=v_{1,3}$ and   by abuse of notation we
also denote by $v$ the image of $v$ in $K(j)$ for $0\le j\le 3$.
Let
$$\overline{K(1)} = K/\bu(\lie g_0[t])\boy_{\boe_1}v =
K(1)/\bu(\lie g_0)\boy_{\boe_1}v.$$

For any $\varepsilon\in\{0,1\}$ and $\bop\in\bz_+^4$
with $\sum_{i=1}^4 p_i\le m_0$, where
$m=3m_0+m_1$ with $0\le m_1\le 2$, set
$$K_{m_1,m_0}^\varepsilon
(\bop)=
 K^{\otimes p_4}\otimes K(2)^{\otimes p_3}\otimes
\overline{K(1)}^{\otimes p_2}\otimes K(1)^{\otimes p_1}\otimes
K(0)^{\otimes m_0-\sum_{i=1}^4p_i}\otimes
K(m_1\omega_1)(\varepsilon).$$

Given an equivalence class  $\bar\bor$ such that $\bar\bor\cap \cal A_1\ne\emptyset$ we assume that
$\bor=\bor_{j,k,s}$ and  let $r_1,r_4$ be defined as in section \ref{conj}. Then set
$$\bor_0 = \bor + r_4(3\boe_1-\boe_2-\boe_4) = (r_1+3r_4)\boe_1 + r_2\boe_2 + r_3\boe_3,$$
where $r_2 = k-s$ and $r_3=s$. 
For $0\le n\le \#\ \bar\bor\cap\cal A_1$, set
$$\bor_n=\bor_{n-1}-(3\boe_1-\boe_2-\boe_4)$$
and define $\bop_n(\bar\bor)\in\bz_+^4$, $\varepsilon(\bar\bor)\in\{0,1\}$
by:
\begin{enumerit}
\item if $m_1=2$, then
$$\bop_{n}(\bar\bor)=(\delta_{2,r_1}, r_2+n, r_3+n,r_4-n),\ \
\varepsilon(\bar\bor)=1-\delta_{0,r_1},$$
\item if $m_1=1, r_1=2$ and $r_2+r_3+r_4 =m_0-1$
(in particular $\#\ \bar\bor\cap\cal A_1=1$), set
$$\bop_0(\bar\bor)=(0, r_2,  r_3+1, r_4),\ \
\varepsilon(\bar\bor)=0,$$
\item and in all other cases,
$$\bop_n(\bar\bor)=( r_1, r_2+n, r_3+n, r_4-n),\ \
\varepsilon(\bar\bor)=0.$$
\end{enumerit}
It is now tedious but not hard  to see that the modules
$K_{m_1,m_0}^{\varepsilon(\bar\bor)}(\bop_n(\bar\bor))$ are defined for
all $0\le n\le \#\ \bar\bor\cap\cal A_1$. Finally, let $v_{\bop_n(\bar\bor)}$ be the
image of the tensor product of the elements $v^{\otimes m_0}\otimes v_{1,m_1}$ in
$K_{m_1,m_0}^{\varepsilon(\bar\bor)}(\bop_n(\bar\bor))$.

\begin{prop}\label{tp}  Let $\bor\in\cal A_1$ be as above and consider $K_{m_1,m_0}^{\varepsilon(\bar\bor)}(\bop_n(\bar\bor))$
for $0\le n\le \#\ \bar\bor\cap\cal A_1$. Write $\bop_n(\bar\bor)=(p_1,p_2,p_3,p_4)$ and $p_0=m_0-\sum_{i=1}^4p_i$.
If $m_1=1, r_1=2$ and $r_2+r_3+r_4 =m_0-1$ we have
\begin{align}\label{crux1gii}
\boy_{\bor_0}v_{\bop_0(\bar\bor)}   =
\boy_{3r_4\boe_1}v^{\otimes r_4} \otimes
x^-_{\alpha_1+\alpha_2}\boy_{(r_3+1)\boe_3}v^{\otimes
r_3+1}\otimes \boy_{r_2\boe_2}v^{\otimes r_2}\otimes v_{1,1},
\end{align}
and there exists $a,b\in\bc$ such that
\begin{equation}\label{crux2gii}
\lie n_0^+ \left(\boy_{\bor_0} +(a x_{\alpha_2}^-x_{\alpha_1}^- +
bx^-_{\alpha_1+\alpha_2})\boy_{\bor_0-2\boe_1+\boe_3}\right)
v_{\bop_0(\bar\bor)} = 0.
\end{equation}
In all other cases we have
\begin{equation}\label{crux0g}
\boy_{\bor_\ell} v_{\bop_n(\bar\bor)}=0 \qquad\text{if}\qquad \ell<n,
\end{equation}
\begin{align}\label{crux1g}
\boy_{\bor_n}v_{\bop_n(\bar\bor)} = \boy_{3p_4\boe_1}v^{\otimes
p_4}\otimes (x_{\alpha_2}^-)^{n}\boy_{p_3\boe_3} v^{\otimes
p_3}\otimes \boy_{p_2\boe_2}v^{\otimes
p_2}\otimes\boy_{p_1\boe_1}v^{\otimes p_1} \otimes v^{\otimes
p_0}\otimes \boy_{\varepsilon(\bar\bor)\boe_1}v_{1,m_1},
\end{align}
and there exists $c_1,\cdots, c_n\in\bc^*$   such that
\begin{equation}\label{crux2g}\lie n_0^+\left(\boy_{\bor_n}+\sum_{\ell=1}^n
c_\ell(x_{\alpha_2}^-)^{\ell}\boy_{\bor_n+\ell(\boe_3-\boe_4)}\right)
v_{\bop_n(\bar\bor)}=0.
\end{equation}
\end{prop}

\begin{pf} A straightforward computation using Lemma \ref{fundvan} prove equations \eqref{crux1gii}, \eqref{crux0g},
and \eqref{crux1g}.
To prove \eqref{crux2g}, let $\bor_{n,\ell} =
\bor_n+\ell(\boe_3-\boe_4)$ for $0\le \ell\le n$. Lemma
\ref{fundvan} now gives,
\begin{align*}
\boy_{\bor_{n,\ell}}v_{\bop_n(\bar\bor)}  =
\boy_{3p_4\boe_1}v^{\otimes p_4}\otimes
(x_{\alpha_2}^-)^{n-\ell}\boy_{p_3\boe_3} v^{\otimes p_3}\otimes
\boy_{p_2\boe_2}v^{\otimes
p_2}\otimes\boy_{p_1\boe_1}v_{1,3}^{\otimes p_1} \otimes
v^{\otimes p_0}\otimes \boy_{\varepsilon(\bar\bor)\boe_1}v_{1,m_1},
\end{align*}
and also that
$$x_{\alpha_1}^+\boy_{\bor_{n,\ell}}v_{\bop_n(\bar\bor)}=\ \
x_{\alpha_2}^+\boy_{3p_4\boe_1}v^{\otimes p_4} =
x_{\alpha_2}^+\boy_{p_1\boe_1}v^{\otimes p_1} =
x_{\alpha_2}^+v^{\otimes p_0} =
x_{\alpha_2}^+\boy_{\varepsilon(\bar\bor)\boe_1}v_{1,m_1} = 0.$$
Since
$$\wt(\boy_{3p_4\boe_1}v^{\otimes p_4})(h_{\alpha_2}) =
\wt(\boy_{p_1\boe_1}v^{\otimes p_1})(h_{\alpha_2}) =
\wt(v^{\otimes p_0})(h_{\alpha_2}) =
\wt(\boy_{\varepsilon(\bar\bor)\boe_1}v_{1,m_1})(h_{\alpha_2}) =0,$$
it follows now that to prove \eqref{crux2g}, it suffices to find $c_1,\cdots
,c_n\in\bc^\times$ such that
$$x^+_{\alpha_2}\left((x_{\alpha_2}^-)^{n}\boy_{p_3\boe_3}v^{p_3} \otimes
\boy_{p_2\boe_2}v^{\otimes p_2} +\sum_{\ell=1}^n c_\ell
(x_{\alpha_2}^-)^{\ell}\left((x_{\alpha_2}^-)^{n-\ell}\boy_{p_3\boe_3}v^{p_3}
\otimes \boy_{p_2\boe_2}v^{\otimes p_2}\right)\right)=0$$
in $K(2)^{\otimes p_3}\otimes \overline{K(1)}^{\otimes p_2}$.
Since  $\wt(\boy_{\boe_3}v) = \wt(\boy_{\boe_2}v) = \omega_2,$ $x_{\alpha_2}^+\boy_{\boe_3}v = 0$, and
in $\overline K(1)$ we have $\boy_{\boe_2}v\ne 0$ and $x_{\alpha_2}^+\boy_{\boe_2}v = 0$,
the result now follows from Lemma \ref{sl2}.

To prove \eqref{crux2gii} we first observe that Lemma \ref{fundvan} also gives
\begin{align*}
\boy_{\bor_0-2\boe_1+\boe_3} v_{\bop_0(\bar\bor)}  =
\boy_{3r_4\boe_1}v_{1,3}^{\otimes r_4} \otimes
\boy_{(r_3+1)\boe_3}v_{1,3}^{\otimes r_3+1}\otimes
\boy_{r_2\boe_2}v_{1,3}^{\otimes r_2}\otimes v_{1,1}.
\end{align*}
The rest of the proof is now similar to the previous case using
Lemma \ref{hwvg2}.
\end{pf}

\subsection{} Let $V$ be any $\lie g_0[t]$--module quotient of $KR(m\omega_1)$ and $v$ be the image of $v_{1,m}$.
Given $\bor\in\bz_+^4$, let
$V^{>\bor}$ be the $\lie g_0$--submodule of $V$
 generated by the elements
$\{\boy_\bos v:\wt(\bos)>\wt(\bor)\}$ and let $V^{\geq\bor}$ be
the $\lie g_ 0$--submodule generated by $V^{>\bor}$ and the
elements $\{\boy_\bos v: \bos\in\bar\bor\}$. For $\mu\in P^+$, and
any finite--dimensional $\lie g_0$--module $W$, let $m_\mu(W)$ be the multiplicity of
the isotypical component in $W$ corresponding to $\mu$.

\begin{prop}\label{ubos} Let $\bos\in\bz_+^4$.  We have,
\begin{eqnarray*} &m_\mu(V^{>\bos})\ne 0 &\implies \mu>\wt(\bos),\\
&m_\mu(V^{\ge \bos})\ne 0 &\implies \mu\ge\wt(\bos).
\end{eqnarray*} In particular,
$$m_\mu(V^{\ge \bos}/V^{>\bos})\ne 0\implies \mu=\wt(\bos).$$

\end{prop}
\begin{pf} Suppose that $m_\mu(V^{>\bos})\ne 0$ for some $\mu\in P^+$.
Let $p_\mu: V^{>\bos}\to V(\mu)^{\oplus m_\mu(V^{>\bos})}$ be the
projection of $\lie g_0$--modules onto the corresponding
isotypical component. Since $p_\mu\ne 0$ it follows that there
must exist $\bor\in\cal A_1$ with $\wt(\bor)>\wt(\bos)$ such that
$p_\mu(\boy_\bor v)\ne 0$. This implies that $\mu-\wt(\bor)\in
Q^+$, i.e., $\mu\ge\wt(\bor)$. The first implication of the Lemma
follows. The second is proved similarly. If
$m_\mu(V^{\ge\bos}/V^{>\bos})\ne 0$, then $p_\mu(\boy_\bor v)\ne
0$ for some $\bor\in \bar\bos$ and hence $\mu\ge \wt(\bos)$. If
$\mu>\wt(\bos)$ then
 there must exist $\bor'$ with $\wt(\bor')=\mu$ such that
$\boy_{\bor'} v$ has non--zero projection onto $V(\mu)$ (see
section \ref{ubpf}). But this is impossible since $\boy_{\bor'}
v\in V^{\bos}$.
\end{pf}

\subsection{Completion of the proof of Proposition \ref{lb} for
$KR(m\omega_1)$}

Given $\bor\in\cal A_1$, let $\bor_n$ and $v_{\bop_n(\bar\bor)}$ be defined as in section \ref{tp}.

\begin{prop}\label{lbpf} Let $\bor\in\cal A_1$ and $u=\sum_{n=0}^{\#\ \bar\bor\cap\cal A_1}c_n\boy_{\bor_n}$ for some $c_n\in\bc$. Then $uv_{1,m}\in  KR(m\omega_1)^{>\bor}$ only if $c_n=0$  for all $0\le n\le \#\ \bar\bor\cap\cal A_1$. In particular, we have
$$KR(m\omega_1)^{\ge\bor}/KR(m\omega_1)^{>\bor}\cong V(\wt(\bor))^{\oplus \ell}$$ for some $\ell\ge \#\ \bar\bor\cap\cal A_1.$
\end{prop}

\begin{pf} Suppose that $c_n\ne 0$ for some $0\le
n\le \#\ \bar\bor\cap\cal A_1$ and assume that $n$ is maximal with this
property. It is not difficult to see that $V = \bu(\lie g_0[t])v_{\bop_n(\bar\bor)}$ is a quotient of $KR(m\omega_1)$.
If $uv_{1,m}\in KR(m\omega_1)^{>\bor}$, then we also have $u v_{\bop_n(\bar\bor)}\in V^{>\bor}$.
On the other hand, it follows from Proposition \ref{tp}
that $u v_{\bop_n(\bar\bor)} = c_n\boy_{\bor_n} v_{\bop_n(\bar\bor)}\ne 0$. But then equations
 \eqref{crux2gii} and \eqref{crux2g}  contradict Lemma \ref{ubos} since $\wt(u v_{\bop_n(\bar\bor)}) =
 \wt(\bor)$. Hence, we must have that $c_n=0$ for all $n$.
\end{pf}

\subsection{The modules $KR^\sigma(m\omega_1)$.}
In what follows we
shall write an element of $\lie g$ as a triple $(x_0,x_1,x_2)$ with
$x_j\in\lie g_j$, $j=0,1,2$.  Let $<,>$ be the Killing form of $\lie
g$. It is not hard to check that the following formulas define an
action of $\lie g[t]^\sigma$ on $K = \lie g_2\oplus\bc$:
\begin{eqnarray*}
&y_0\otimes t^{3r}(x_2,a)&= \delta_{r,0}([y_0,x_2], 0),\\
&y_1\otimes t^{3r+1}(x_2,a)&= \delta_{r,0}(0, <x_2,y_1>a),\\
&y_2\otimes t^{3r+2}(x_2,a)&= 0,
\end{eqnarray*}
where $y_j\in\lie g_j$ for $j=0,1,2$.
Moreover, since $<,>$ is non degenerate on $\lie g_1\times \lie g_2$, it is not hard to see that the assignment
 $z^+_{2\alpha_1+\alpha_2},\mapsto v_{1,1}^\sigma$ extends to an isomorphism of $\lie g[t]^\sigma$--modules
 $K\cong KR^\sigma(\omega_1)$. For $m>1$ we proceed exactly as in section \ref{om2}.

\subsection{The modules $KR^\sigma(m\omega_2)$.} Proceeding as in the previous section we
see that the following formulas define an
action of $\lie g[t]^\sigma$ on $\lie g\oplus\bc$:
\begin{eqnarray*}&y_0\otimes t^{3r}(x_0,x_1,x_2,a)&=
\delta_{r,0}([y_0,x_0],[y_0,x_1], [y_0,x_2], 0)+
(0,0,0,\delta_{r,1}<x_0,y_0>a),\\
&y_1\otimes t^{3r+1}(x_0,x_1,x_2,a)&= \delta_{r,0}(0,[y_1,x_0],
[y_1,x_1], <x_2,y_1>a),\\
&y_2\otimes t^{3r+2}(x_0,x_1,x_2,a)&= \delta_{r,0}(0,0, [y_2,x_0],
<x_1,y_2>a),\end{eqnarray*} where $x_j,y_j\in\lie g_j$ for $j=0,1,2$.
Moreover, it is straightforward to check that this module is a
quotient of $KR^\sigma(\omega_2)$ and hence proves Proposition
\ref{lb} when $m=1$.

For $m>1$ the proof follows the same pattern  as that for $KR(m\omega_1)$,
$m>3$, and we just list the relevant modifications and omit the
details. Let $\bor\in\bz_+^4$. Similarly to Lemma \ref{fundvan} we see that
in $KR^\sigma(\omega_2)$ we have
$$\boy^\sigma_\bor v_{2,1}^\sigma=0\ \ \iff\ \  \bor\notin\{\mathbf 0,\boe_1,\boe_2,\boe_3,\boe_4,\boe_1+\boe_2\},$$

$$\lie n_0^+ \boy^\sigma_{\boe_1}v^\sigma_{2,1}=\lie n_0^+ \boy^\sigma_{\boe_2}v^\sigma_{2,1}=\lie n_0^+ \boy^\sigma_{\boe_4}v^\sigma_{2,1}=0,
\qquad\text{and}\qquad \boy^\sigma_{\boe_3}v_{2,1}^\sigma \in\bc x^-_{\alpha_1}\boy^\sigma_{\boe_2}v^\sigma_{2,1}.$$

Lemma \ref{ubos} is still valid with the obvious modifications.
Set $K = KR^\sigma(\omega_2)$, $v=v^\sigma_{2,1}$ and   by abuse of notation we
also denote by $v$ the image of $v$ in $K(j)$ for $0\le j\le 3$.
For any $\bop=(p_1,p_2,p_3)\in\bz_+^3$ satisfying $p_1+p_2+p_3\le m$ define
$$K_m(\bop) = K^{\otimes p_3}\otimes K(2)^{\otimes p_2}\otimes K(1)^{\otimes p_1}\otimes K(0)^{\otimes m-p_1-p_2-p_3},$$
and let $v_{\bop}$ be the image of $v^{\otimes m}$ in $K_m(\bop)$.

Now fix $\bar\bor$ such that $\bar\bor\cap\cal A_2^\sigma\ne \emptyset$ and assume that $\bor=\bor_{j,k,s} = (r_1,r_2,0,r_4)$, where $\bor_{j,k,s}$ is defined as in section \ref{conj}.
Then, for $0\le n\le \#\ \bar\bor\cap \cal A_2^\sigma$,  set $\bor_n = \bor+n(\boe_1+\boe_3-\boe_4)$ and
$\bop_n(\bar\bor) = (r_1+n,r_2+n,r_4-n)$.
  We have the following analog of Proposition \ref{tp}.

\begin{prop}\label{tpd4} Let $\bor\in\cal A_2^\sigma$ be as above and consider $K_m(\bop_n(\bar\bor))$ for
$0\le n\le \#\ \bar\bor\cap \cal A_2^\sigma$.
 Then
$\boy^\sigma_{\bor_\ell}v_{\bop_n(\bar\bor)}=0$ if $ \ell<n$,
\begin{equation*}\label{crux1}
\boy^\sigma_{\bor_n}v_{\bop_n(\bar\bor)}= \boy^{\sigma}_{(r_4-n)\boe_4}v^{\otimes r_4-n} \otimes
(x^-_{\alpha_1})^n\boy^\sigma_{(r_2+n)\boe_2}v^{\otimes r_2+n}
\otimes \boy^\sigma_{(r_1+n)\boe_1} v^{\otimes r_1+n}\otimes
v^{\otimes m-r_4-r_2-r_1-n},
\end{equation*}
and there exist  $c_1,\cdots, c_n\in\bc^\times$ such that
\begin{equation*}\label{crux2}
\lie n_0^+\left(\boy^\sigma_{\bor_n}+\sum_{\ell=1}^n
c_\ell(x_{\alpha_1}^-)^\ell\boy^\sigma_{\bor_n+\ell(\boe_2-\boe_3)}\right)v_{\bop_n(\bar\bor)}=0.
\end{equation*}
\hfill\qedsymbol
\end{prop}

The proof of Proposition \ref{lb} is then completed as
before by using the obvious modification of Proposition
\ref{lbpf}.

\bibliographystyle{amsplain}

\end{document}